\documentclass [12pt]{article}

\title {A dichotomy for $D$-rank 1 types in simple theories}
\author {Ziv Shami}
\newtheorem {theorem}{Theorem}[section]
\newtheorem {lemma}[theorem]{Lemma}
\newtheorem {definition}[theorem]{Definition}

\newtheorem {fact}[theorem]{Fact}

\newtheorem {corollary}[theorem]{Corollary}

\newtheorem {remark}[theorem]{Remark}

\newtheorem {proposition}[theorem]{Proposition}

\newtheorem {notation}[theorem]{Notation}
\newtheorem {claim}[theorem]{Claim}

\def\proof {\noindent \textbf{Proof:} }
\def\qed {$\ \ \ \ \Box$\\}

\def\qedend {$\ \ \ \ \ \ \ \ \ \ \ \ \ \ \ \ \ \ \ \ \ \ \ \ \ \ \ \ \ \ \ \ \ \ \ \ \ \ \ \ \ \ \ \ \ \ \ \ \ \ \
 \ \ \ \ \ \ \ \ \ \ \ \ \ \ \ \ \ \ \ \ \ \ \ \ \ \ \ \ \ \ \ \ \ \ \ \ \ \ \ \ \ \ \ \Box$\\}

\newsavebox{\indbin}
\savebox{\indbin}{\begin{picture}(0,0)
\newlength{\gnu}
\settowidth{\gnu}{$\smile$} \setlength{\unitlength}{.5\gnu} \put(-1,-.65){$\smile$}
\put(-.25,.1){$|$}
\end{picture}}
\newcommand{\nonfork}[3]
{\mbox{$\begin{array}{ccc} \mbox{$#1$} & \usebox{\indbin} & \mbox{$#2$} \\
        & \mbox{$#3$} &
\end{array}$}}
\newcommand{\nonforkempty}[2]
{\mbox{$\begin{array}{ccc} \mbox{$#1$} & \usebox{\indbin} & \mbox{$#2$}
\end{array}$}}
\newcommand{\fork}[3]
{\mbox{$\begin{array}{ccc} \mbox{$#1$} & \!\mbox{$\!\!\not\!\:\usebox{\indbin}$} & \mbox{$#2$} \\
        & \mbox{$#3$} &
\end{array}$}}

\newsavebox{\sindbin}
\savebox{\sindbin}{\begin{picture}(0,0)
\newlength{\sgnu}
\settowidth{\sgnu}{$\smile$} \setlength{\unitlength}{.5\sgnu} \put(-1,-.65){$\smile$}
\put(-.25,.1){$|s$}
\end{picture}}

\newsavebox{\starindbin}
\savebox{\starindbin}{\begin{picture}(0,0)
\newlength{\stargnu}
\settowidth{\stargnu}{$\smile$} \setlength{\unitlength}{.5\stargnu} \put(-1,-.65){$\smile$}
\put(-.25,.1){$|*$}
\end{picture}}

\newsavebox{\qindbin}
\savebox{\qindbin}{\begin{picture}(0,0)
\newlength{\qgnu}
\settowidth{\qgnu}{$\smile$} \setlength{\unitlength}{.5\qgnu} \put(-1,-.65){$\smile$}
\put(-.25,.1){$|_{qf}$}
\end{picture}}

\def\card #1 {{\vert #1 \vert}}

\def\CC {{\cal C}}

\def\PP {{\cal P}}
\def\UU {{\cal U}}

\usepackage{amsfonts}

\begin{document}

\maketitle

\begin{abstract}
We prove a dichotomy for $D$-rank 1 types in simple theories that generalizes Buechler's dichotomy
for $D$-rank 1 minimal types in stable theories: every $D$-rank 1 type is either 1-based or part of
its algebraic closure, defined by a single formula, almost contains a non-algebraic formula that
belongs to a non-forking extension of the type. In addition we prove that a densely 1-based type of
$D$-rank 1 is 1-based. We also observe that for a hypersimple unidimensional theory the existence
of a non-algebraic stable type implies stability (and thus superstability).
\end{abstract}

\section{Introduction}
In 1985 Buechler proved [B] a remarkable dichotomy between model theoretic simplicity and geometric
simplicity; it says that any minimal $D$-rank 1 type in a stable theory is either 1-based or has
Morley rank 1. In this paper we give a generalization of this result for any $D$-rank 1 type of an
arbitrary simple theory. As a special case we get Buechler's dichotomy for any $D$-rank 1 minimal
type in exactly the form mentioned above. The proof applies certain properties of the forking
topology introduced in [S1], a variant of the topologies introduced in [H0,P]. In these papers,
this topology (and generalizations of it in [S1]) has been used to obtain certain approximations of
definable sets of finite rank for proving supersimplicity of countable hypersimple/hypersimple
low/stable unidimensional theories. It is the minimal topology on $S_x(A)$ such that all the
relations $\Gamma_F(x)$ defined by $\Gamma_F(x)=\exists y (F(x,y)\wedge \nonfork{y}{x}{A})$ are
closed for any type-definable relation $F(x,y)$ over $A$. In the current paper we show in addition
that the notion of an essentially 1-based type introduced in [S1] coincides with the notion of a
1-based type in the $D$-rank 1 case. In fact, we show this for a more general class of $D$-rank 1
types that we call densely 1-based types. A corollary of this is that in any supersimple
unidimensional theory (e.g. any countable hypersimple unidimensional theory [S1]), if for any small
set $A$, finite tuple $c$ and sort $s$, the set $\{a\in \CC^s \vert\nonfork{a}{A\bar c}{acl(aA)\cap
acl(\bar cA)} \}$ is dense in the forking topology on $S_s(A\bar c)$ then $T$ is 1-based. A
posteriori, our result shows that the case handled in [S1] in which the theory is essentially
1-based is in fact just the case in which it is 1-based (for which the proof is much easier) but,
of course, that doesn't mean we can simplify the proof as we don't know a $D$-rank 1 type exists.

The notations are standard, and throughout the paper we work in a highly saturated, highly strongly
homogeneous model $\CC$ of a complete first-order theory $T$ in a language $L$ with no finite
models. We will often work in $\CC^{eq}$.

\section{Preliminaries}
We assume basic knowledge of simple theories as in [K],[KP],[HKP] as well as some knowledge on
hyperimaginaries, (almost-) internality and analyzability in simple theories (we follow the
terminology in [W, 3.4]). In this section, we recall some basic facts related to the forking
topology and to pairs of models in a simple theory that are relevant for this paper. In this
section $T$ will denote a simple theory and we work in $\CC$.

\subsection {The forking topology}

\begin{definition} \label{tau definition}
\em Let $A\subseteq \CC$ and let $x$ be a finite tuple of variables. An invariant set $\UU$ over
$A$ is said to be \em a basic $\tau^f$-open set over $A$ \em if there is a $\phi(x,y)\in L(A)$ such
that $$\UU=\{a \vert \phi(a,y)\ \mbox{forks\ over}\ A \}.$$
\end{definition}

\noindent Note that the family of basic $\tau^f$-open sets over $A$ is closed under finite
intersections, thus form a basis for a unique topology on $S_x(A)$ which we call the
$\tau^f$-topology or the forking-topology.

\begin{remark} \label{ft_remark}\em
Note that the forking-topology on $S_x(A)$ refines the Stone-topology (for every $x$ and $A$) and
that $\{a\in\CC^x\vert a\not\in acl(A)\}(=\{a\in\CC^x\vert x=a\ \mbox{forks\ over}\ A\})$ is a
forking-open subset of $S_x(A)$ (when we identify $A$-invariant sets with subsets of $S_x(A)$).
\end{remark}

\begin{definition}\label{projection closed}
\em We say that \em the $\tau^f$-topologies over $A$ are closed under projections ($T$ is PCFT over
$A$) \em if for every $\tau^f$-open set $\UU(x,y)$ over $A$ the set $\exists y \UU(x,y)$ is a
$\tau^f$-open set over $A$. We say that \em the $\tau^f$-topologies are closed under projections
($T$ is PCFT) \em if they are over every set $A$.
\end{definition}

\begin{fact}\label{tau extensions}[S0, Lemma 2.6]
Let $\UU$ be a $\tau^f$-open set over a set $A$ and let $B\supseteq A$ be any set. Then $\UU$ is
$\tau^f$-open over $B$.
\end{fact}

We say that an $A$-invariant set $\UU$ has $SU$-rank $\alpha$ and write $SU(\UU)=\alpha$ if
$Max\{SU(p) \vert p\in S(A), p^\CC\subseteq\UU\}=\alpha$. We say that $\UU$ \em has bounded finite
$SU$-rank if $SU(\UU)=n$ for some $n<\omega$. \em

\begin{fact}\label{tau bounded SU}[S0, Lemma 2.12, Proposition 2.13]
Let $\UU$ be an unbounded $\tau^f$-open set over some set $A$. Assume $\UU$ has bounded finite
$SU$-rank. Then there exists a set $B\supseteq A$ and $\theta(x)\in L(B)$ of $SU$-rank 1 such that
$\theta^\CC\subseteq \UU\cup acl(B)$. In case $SU(\UU)=1$, the set $acl_s(A)\cup \UU$ is
Stone-open, where $s$ is the sort of (elements of) $\UU$.
\end{fact}

\subsection {The extension property being first-order and PCFT}

We recall some natural extensions of notions from [BPV]. By a pair $(M,P^M)$ of $T$ we mean an
$L_P=L\cup \{P\}$-structure, where $M$ is a model of $T$ and $P$ is a new predicate symbol whose
interpretation is an elementary submodel of $M$.  For the rest of this subsection, by a $\vert T
\vert$-\em small type \em we mean a complete hyperimaginary type in $\leq \vert T\vert$ variables
over a hyperimaginary of length $\leq \vert T\vert$ (i.e. a sequence of length $\leq \vert T\vert$
modulo a $\emptyset$-type-definable equivalence relation).

\begin{definition}\label {def1}\em
Let $\PP_0,\PP_1$ be $\emptyset$-invariant families of $\vert T\vert$-small types.

\noindent 1) We say that a pair $(M,P^M)$ satisfies the extension property for $\PP_0$ if for every
$L$-type $p\in S(A)$, where $A$ is a hyperimaginary with $A\in dcl(M)$ and $p\in \PP_0$, there is
$a\in p^M$ such that $\nonfork{a}{P^M}{A}$.

\noindent 2) Let $$T_{Ext,\PP_0}=\bigcap \{Th_{L_P}(M,P^M) \vert\ \mbox{the\ pair\ $(M,P^M)$
satisfies\ the\ extension\ property\  w.r.t.\  } \PP_0\ \}.$$

\noindent 3) We say that $\PP_0$ dominates $\PP_1$ w.r.t. the extension property if $(M,P^M)$
satisfies the extension property for $\PP_1$ for every $\vert T\vert^+$-saturated pair
$(M,P^M)\models T_{Ext,\PP_0}$. In this case we write $\PP_0\unrhd_{_{Ext}} \PP_1$.

\noindent 4) We say that the extension property is first-order for $\PP_0$ if
$\PP_0\unrhd_{_{Ext}}\PP_0$ (i.e. every $\vert T\vert^+$-saturated model of $T_{Ext,\PP_0}$
satisfies the extension property for $\PP_0$). We say that the extension property is first-order if
the extension property is first-order for the family of all $\vert T\vert$-small types
(equivalently, for the family of all real types over sets of size $\leq \vert T\vert$).
\end{definition}

\begin{fact} \label{ext family}[S1, Lemma 3.7]
Let $\PP_0$ be an $\emptyset$-invariant family of $\vert T\vert$-small types. Assume $\PP_0$ is
extension-closed and that the extension property is first-order for $\PP_0$. Let $\PP^*$ be the
maximal class of $\vert T\vert$-small types such that $\PP_0\unrhd_{_{Ext}}\PP^*$. Then
$\PP^*\supseteq An(\PP_0)$, where $An(\PP_0)$ denotes the class of all $\vert T\vert$-small types
analyzable in $\PP_0$ by hyperimaginaries.
\end{fact}

\begin{fact}\label{ext pcft}[S1, Corollary 3.13]
Suppose the extension property is first-order in $T$. Then $T$ is PCFT.
\end{fact}

\section{Dichotomies for rank 1 types}
We first prove a dichotomy between essential 1-basedness and strong minimality for any minimal type
$p$ with possibly no ordinal $D$-rank. In fact, this is a special case: if merely $p$ has $SU$-rank
1, there still is a dichotomy but we don't get a strongly minimal set. Then we prove a strong
version of this for $D$-rank 1 types:  any such type is either 1-based or its algebraic closure (by
a single formula) almost contains a non-algebraic formula; in the special case when $p$ is in
addition stationary (thus minimal) we conclude that if $p$ is not 1-based then it has Morley rank
1. In this section $T$ is assumed to be a simple theory and we work in $\CC^{eq}$ (in particular,
all types are assumed to be
types of imaginaries unless otherwise stated.\\

The following definability result [S1, Proposition 4.4] will be useful.

\begin{fact}\label{Open_Cb}
Let $q(x,y)\in S(\emptyset)$ and let $\chi(x,y,z)\in L$ be such that $\models \forall y\forall
z\exists^{<\infty} x\chi(x,y,z)$. Then the set $$\UU=\{(e,c,b,a) \vert\ e\in acl(Cb(cb/a))\}$$ is
relatively Stone-open inside the type-definable set
$$F=\{(e,c,b,a)\vert\ \nonforkempty{b}{a}, \models\chi(c,b,a), tp(cb)=q\}.$$ (where
$e$ is taken from a fixed sort too).
\end{fact}

\begin{remark}\label{Open_Cb_remark}\em
Note that Fact \ref{Open_Cb} is equivalent to the assertion that the fiber of $\UU$ obtained by
fixing a tuple $cb$ (in the second and third coordinates of $\UU$) whose type is $q$ is relatively
Stone-open inside the corresponding fiber of $F$. This is true simply because $q$ is complete (Fact
\ref{Open_Cb} trivially implies that each such fiber is relatively Stone-open).
\end{remark}

First, we prove a certain version of the dichotomy theorem from [S1] (and generalizations of it in
[S2]) that is closely related to Buechler's dichotomy. Here we assume that $p$ itself is not
essentially 1-based (rather than some type that is internal in $p$) and find a non-algebraic
definable set contained in the algebraic closure of $p$ (rather than almost-internal in $p$).
Recall that a partial type $p$ is called minimal if for every set $A$ containing its domain there
is a unique complete non-algebraic type over $A$ extending $p$. A weakly minimal definable set is
one defined by a $D$-rank 1 formula. We recall now the definition of an essentially 1-based type
[S1]:

\begin{definition}\label {def ess-1-based}\em
1) A type $p\in S(A)$ is said to be \em essentially 1-based \em if for every finite tuple $\bar c$
from $p$ and for every type-definable forking-open set $\UU$ over $A\bar c$, the set $\{a\in \UU
\vert\ Cb(a/A\bar c)\not\in bdd(aA)\}$ is nowhere dense in the Stone-topology of $\UU$
\footnote{Note that if $\UU=\emptyset$, the requirement on $\UU$ is trivially satisfied.}.
\end{definition}

Recall that a theory $T$ is called \em hypersimple \em if it is simple and eliminates
hyperimaginaries [S1].

\begin{proposition}
Let $T$ be a countable hypersimple theory and assume $T^{eq}$ has PCFT. Let $p\in S(\emptyset)$ be
a type of $SU$-rank 1 that is not essentially 1-based by means of the forking-topology. Then
$acl(p^\CC)$ contains a weakly minimal definable set defined over $acl(p^\CC)$. If, in addition,
$p$ is minimal then $acl(p^\CC)$ contains a strongly minimal definable set.
\end{proposition}

\proof By the assumption, there exists a finite tuple $\bar c$ of realizations of $p$, a
type-definable forking-open set $\UU$ over $\bar c$ and a $\bar c$-invariant Stone-dense set
$D\subseteq \UU$ such that $Cb(\bar c/a)\not\subseteq acl(\bar c)$ for all $a\in D$ (note that
$Cb(\bar c/a)\subseteq acl(\bar c)$ iff $\nonfork{\bar c}{a}{acl(\bar c)\cap acl(a)}$ iff
$Cb(a/\bar c)\subseteq acl(a)$).

\begin{claim}\label{claim_baire}
There are disjoint tuples $\bar c_0,\bar c_1$ such that $\bar c=\bar c_0\cup \bar c_1$ and
$\chi(\bar x_1,\bar x_0,y)\in L$ with $\forall \bar x_0 y \exists^{<\infty}\bar x_1\chi(\bar
x_1,\bar x_0,y)$ such that

$$\UU'=\{a\in \UU \vert \nonforkempty{a}{\bar c_0},\ \chi(\bar c_1,\bar c_0,a)\}$$ contains a non-empty
invariant set over $\bar c$ that is relatively Stone-open in $\UU$.
\end{claim}

\proof Let $S$ be the set of triples $(\bar c_0,\bar c_1,\chi)$ such that $\bar c=\bar c_0\cup \bar
c_1$, $\ \bar c_0,\bar c_1$ are disjoint and $\chi=\chi(\bar x_1,\bar x_0,y)\in L$ is such that
$\forall \bar x_0 y \exists^{<\infty}\bar x_1\chi(\bar x_1,\bar x_0,y)$. For every $(\bar c_0,\bar
c_1,\chi)\in S$, let
$$F_{\bar c_0,\bar c_1,\chi}=\{a\in \UU\vert\ \nonforkempty{a}{\bar c_0},\ \chi(\bar c_1,\bar c_0,a)\}.$$ Since $SU(p)=1$, we
get that $\UU=\bigcup_{_{(\bar c_0,\bar c_1,\chi)\in S}} F_{\bar c_0,\bar c_1,\chi}$. Since each
$F_{\bar c_0,\bar c_1,\chi}$ is Stone-closed in $\UU$ and $S$ is countable (as $L$ is countable),
the claim follows by the Baire category theorem for the Stone-topology of $\UU$.$\ \ \ \ \Box$\\

\noindent By Claim \ref{claim_baire}, there is $\UU_0\subseteq \UU'$ that is non-empty, Stone-open
in $\UU$, type-definable and forking-open over $\bar c$ (intersect $\UU$ with an appropriate
definable set over $\bar c$). Since $D$ is Stone-dense in $\UU$, there exists $a^*\in D\cap\UU_0$.
Then $\hat e=Cb(\bar c/a^*)\not\subseteq acl(\bar c)$. Let $e^*\in dcl(\hat e)\backslash acl(\bar
c)$. Let $s^*$ be the sort of $e^*$.

\begin{claim}\label{claim_def}
The set $\{(e,a)\ \vert\ e\in acl(Cb(\bar c/a))\}$ is relatively Stone-open over $\bar c$ inside
$\CC^{s^*}\times \UU_0$.
\end{claim}

\proof This is a conclusion of Remark \ref{Open_Cb_remark} and the fact that $\UU_0\subseteq \UU'$.
$\ \ \ \ \Box$

\begin{claim}\label{claim_def1}
Let $E=\{e\in \CC^{s^*}\vert\ \exists a\in\UU_0 [e\in acl(Cb(\bar c/a))]\}.$ There exists an
unbounded forking-open type-definable set $E'\subseteq E$ over $\bar c$.
\end{claim}

\proof By Claim \ref{claim_def}, there exists a Stone-open set $V=V(x,y)$ over $\bar c$ such that
$V\cap (\CC^{s^*}\times \UU_0)=\{(e,a)\in \CC^{s^*}\times \UU_0 \vert\ e\in acl(Cb(\bar c/a))\}$.
Clearly, there exists a definable $V_0\subseteq V$ such that $(e^*,a^*)\in V_0\cap (\CC^{s^*}\times
\UU_0)$. Since $T^{eq}$ has PCFT, $E'=\{e\in \CC^{s^*}\vert\ \exists a\in\UU_0(V_0(e,a))\}$ is an
unbounded forking-open set over $\bar c$. Since $\UU_0,V_0$ are type-definable , $E'$ is
type-definable. Clearly, $E'\subseteq E$.$\ \ \ \ \Box$\\

Since every element of $E'$ is in the algebraic closure of some finite tuple of realizations of
$p$, we may assume, by the Baire category theorem for the Stone topology of $E'$, that $SU(E')=n^*$
for some $0<n^*<\omega$. Let $e'\in E'$ be such that $SU(e'/\bar c)=n^*$ and let $A\supseteq\bar c$
be such that $SU(e'/A)=n^*-1$. By passing to the canonical base of $Lstp(e'/A)=stp(e'/A)$ we may
assume that $A\subseteq acl(p^\CC)$ (we may assume $A$ is finite). Let $E'_1=\{e_1'\in E'\vert
\phi(e_1',a)\}$, where $\phi(x,a)\in L(A)$ is a formula in $tp(e'/A)$ that forks over $\bar c$.

Note that since $E'$ is a forking-open set over $\bar c$, it is a forking-open set over $A$ by Fact
\ref{tau extensions} (being forking-open over $A$ is not an immediate corollary of being a
forking-open set over some subset of $A$). Now, $E'_1$ is a forking-open set over $A$ (and
type-definable) as it is the intersection of a forking-open set over $A$ and a definable set over
$A$. Clearly,  $SU(E'_1)=n^*-1$. By repeating this we get a forking-open (type-definable) set
$E^*\subseteq acl(p^\CC)$ over a finite set $A^*\subseteq acl(p^\CC)$ of $SU$-rank 1. By Fact
\ref{tau bounded SU}, there exists a non-algebraic formula $\theta(x)\in L(A^*)$ such that
$\theta(\CC)\subseteq E^* \cup acl(A^*)\subseteq acl(p^\CC)$. Clearly, $\theta(x)$ is weakly
minimal. If $p$ is minimal then $\theta(x)$ has ordinal Morley rank (since the language is
countable, every realization of $\theta(x)$ is in the algebraic closure of some tuple of
realizations of $p$ and thus by minimality of $p$, for every countable set $A$, the number of types
of realizations of $\theta$ over $A$ is countable). Thus there exists a strongly minimal
$\theta^*(x)\vdash \theta(x)$.$\ \ \ \ \Box$

\begin{theorem}\label{main}
Let $p\in S(\emptyset)$ be a type of $D$-rank 1. Then either $p$ is 1-based or there exists
$\tilde\chi(x,\bar z)\in L$ with $\forall \bar z\exists ^{<\infty}x\tilde\chi(x,\bar z)$ and an
$\emptyset$-independent tuple $\bar c$ of realizations of $p$, and a (non-algebraic) formula
$\theta(x)\in L(\bar c)$ in some non-forking extension $\bar p$ of $p$ such that for any
non-algebraic realization $a$ of $\theta(x)$ there is an $\emptyset$-independent tuple $\bar c'$ of
realizations of $p$ such that $\tilde\chi(a,\bar c')$.
\end{theorem}

Before presenting the proof we recall some standard terminology.

\begin{definition}\em
A type $p\in S(\emptyset)$ is \em 1-based \em if for every set $C$ and tuple $\bar a\subseteq
p^\CC$ we have $Cb(\bar a/C)\in bdd(\bar a)$.
\end{definition}

\begin{definition}\em
An $SU$-rank 1 type $p\in S(\emptyset)$ is called \em linear \em if for every set $C$ and all
$a,b\in p^\CC$ with $SU(ab/C)=1$ we have $SU(Cb(ab/C))\leq 1$.
\end{definition}

\begin{fact}\label{linearity}[V,DK]
Assume $p\in S(\emptyset)$ is a type of $SU$-rank 1. Then $p$ is 1-based iff $p$ is linear.
\end{fact}

\begin{lemma}\label{lemma1}
Let $p\in S(\emptyset)$ be a type of $SU$-rank 1. Then $p$ is 1-based iff for every $a,b\in p^\CC$
and finite tuple $\bar c$ of realizations of $p$ we have $\nonfork{ab}{\bar c}{bdd(ab)\cap bdd(\bar
c)}$.

\end{lemma}

\proof First, note the following general observation (an easy $SU$-rank calculation).

\begin{claim}\label{SU cal}
Assume $SU(a)=2$ and $SU(a/C)=1$. Then $\nonfork{a}{C}{bdd(a)\cap bdd(C)}$ iff $SU(Cb(a/C))=1$.
\end{claim}

\noindent Now, clearly we only need to prove right to left. Assume the right hand side holds. By
Claim \ref{SU cal} and Fact \ref{linearity}, it will be sufficient to show that for every $a,b\in
p^\CC$ and set $C$ such that $SU(ab/C)=1$ we have $\nonfork{ab}{C}{bdd(ab)\cap bdd(C)}$ (we may
clearly assume $SU(ab)=2$, i.e. $a$ is independent of $b$ over $\emptyset$, since otherwise
$SU(ab/C)=1$ implies that $ab$ is independent of $C$). Indeed, otherwise $e\equiv Cb(ab/C)\not\in
bdd(ab)$. Let $(a_ib_i \vert i<\omega)$ be a sequence such that $(a_ib_i \vert i<\omega)^\wedge ab$
is a Morley sequence of $Lstp(ab/C)$. Then $e\in dcl(a_ib_i \vert i<\omega)$, and therefore
$\fork{ab}{(a_ib_i \vert i<\omega)}{bdd(ab)\cap bdd(a_ib_i \vert i<\omega)}$. By passing to a
finite subsequence of $(a_ib_i \vert i<\omega)$ we get a contradiction to our assumption.$\ \ \ \ \Box$\\


\begin{remark}\label{remark1}\em
Let $T$ be a 1-sorted theory of $SU$-rank 1. Then the extension property is first order in $T$ and
in $T^{eq}$; thus $T$ and $T^{eq}$ have PCFT.
\end{remark}

\proof By [H1], any 1-sorted theory of $SU$-rank 1 eliminates the $\exists^{\infty}$ quantifier.
Thus the extension property is first-order for 1-types (see [V, Proposition 2.15]). Since every
non-algebraic type is non-orthogonal to a 1-type, the extension property is first-order in $T$ and
in $T^{eq}$ by Fact \ref{ext family}. By Fact \ref{ext pcft}, $T$ and $T^{eq}$ have PCFT.$\ \ \ \ \Box$\\

\begin{notation}\em
Let $D$ be a definable set over $\emptyset$. Let $D_*$ be the induced structure on $D$, that is,
the universe of $D_*$ is $D$ and it is equipped with all $\emptyset$-definable subsets of $\CC$
that are subsets of $D^n$ for some $n<\omega$. $D_*$ is saturated. Note that $D_*^{eq}$ can be
interpreted as the induced structure of $\CC^{eq}$ on $D$ and appropriate disjoint
$\emptyset$-definable subsets of $\CC^{eq}$.
\end{notation}

\begin{remark}\label{claim1}\em
Let $T$ be any simple theory. Let $D$ be a non-algebraic $\emptyset$-definable set. Then for every
tuples $\bar a, \bar b$ from $D$, we have $\CC\models \nonfork{\bar a}{\bar b}{bdd(\bar a)\cap
bdd(\bar b)}$ iff $D_*\models \nonfork{\bar a}{\bar b}{bdd(\bar a)\cap bdd(\bar b)}$.
\end{remark}

\proof Let $e_\CC=Cb^{\CC}(\bar a/\bar b)$ be the canonical base of $Lstp(\bar a/\bar b)$ in the
sense of $\CC$, and let $e_{D_*}=Cb^{D_*}(\bar a/\bar b)$ be the canonical base of $Lstp(\bar
a/\bar b)$ in the sense of $D_*$. We need to show that $\CC\models e_{\CC}\in bdd(\bar a)$ iff
$D_*\models e_{D_*}\in bdd(\bar a)$. Now, note that for every partial type $p(x,c)$ of $D_*$ and
small set $A$ of $D_*$, $p(x,c)$ doesn't fork over $A$ in the sense of $D_*$ iff $p(x,c)$ doesn't
fork over $A$ in the sense of $\CC$. Thus $e_{\CC}=e_{D_*}$ and the claim follows.$\ \ \ \ \Box$\\

\noindent\textbf{Proof of Theorem \ref{main}:} Assume $p$ is not 1-based. Let $D\in p$ be of
$D$-rank 1. By Lemma \ref{lemma1},  there exists $a,b\in p^\CC$ and a finite tuple $\bar
c=c_0c_1...c_n$ of realizations of $p$ such that $\fork{ab}{\bar c}{bdd(ab)\cap bdd(\bar c)}$. By
Remark \ref{claim1}, $D_*\models \fork{ab}{\bar c}{bdd(ab)\cap bdd(\bar c)}$. \textbf{From now on
we work in $D_*$}. As $SU(p)=1$, we may clearly assume $\bar c$ is an $\emptyset$-independent
sequence of realizations of $p$. Clearly, $\nonforkempty{a}{\bar c}$, $\chi_0(b,a,\bar c)$ for some
$\chi_0(x,y,\bar z)\in L$ such that $\forall y\forall \bar z\exists ^{<\infty}x\chi_0(x,y,\bar z)$,
and we may assume $\chi_1(c_n,c_{n-1},...,c_0,a,b)$ for some $\chi_1\in L$ such that
$$\forall z_0z_1...z_{n-1}\forall xy \exists ^{<\infty}z_n\chi_1(z_n,z_{n-1},...z_0,x,y).$$
Therefore $\{c_0,c_1,...,c_{n-1},a,b\}$ is $\emptyset$-independent (as the dimension of
$\{c_0,c_1,c_2,...,c_n,a,b\}$ in the pregeometry $(p^\CC,acl)$ is $n+2$.) Let $$\tilde
D=\{(a',b')\in D^2\vert\ a'\not\in acl(\bar c), \chi_0(b',a',\bar c),
\chi_1(c_n,c_{n-1},...,c_0,a',b')\}.$$ By Remark \ref{ft_remark}, $\tilde D$ is a forking-open set
over $\bar c$, and $(a,b)\in\tilde D$.

\begin{claim}\label{subclaim1}
Let $\alpha$ be any ordinal. There is a partial type $\Lambda(xy,\langle\bar z_i \vert
i<\alpha\rangle)$ over $\bar c$ such that for every $(a',b')\in\tilde D$, for every sequence
$\langle\bar c_i\vert i<\alpha\rangle$, we have $\Lambda(a'b',\langle\bar c_i\vert
i<\alpha\rangle)$ iff $\langle\bar c_i\vert i<\alpha\rangle$ is a Morley sequence of $tp(\bar
c/a'b')$ that starts at $\bar c$.
\end{claim}

\proof By the definition of $\tilde D$, for every $(a',b')\in\tilde D$, a sequence $\langle\bar
c_i\vert i<\alpha\rangle$ that starts at $\bar c$ is a Morley sequence of  $tp(\bar c/a'b')$ iff it
is indiscernible over $a'b'$ and $\{a'b'\}\cup \{\bar c_i^{<n}\vert i<\alpha\}$ is independent over
$\emptyset$. Since the type of $\bar c_i^{<n}$ in such a sequence is fixed (equal to $tp(\bar
c^{<n})$) the required condition is a type-definable condition over $\bar c$ on $(a',b')\in\tilde
D$ and sequence
$\langle\bar c_i\vert i<\alpha\rangle$.$\ \ \ \ \Box$\\

\begin{claim}\label{subclaim2}
There exists a $\chi^*(xy,\bar z)\in L$ and $m^*<\omega$ such that $\forall \bar z \exists
^{<\infty}xy\ \chi^*(xy,\bar z)$ and such that for every Morley sequence $(\bar c_i\vert i\leq
m^*)$ of $tp(\bar c/ab)$ that starts at $\bar c$ we have $$\bigvee_{0\leq i<j\leq m^*}
\chi^*(ab,\bar c_0^{<n},\bar c_1^{<n},..., \bar c_{m^*}^{<n}, c_i^n,c_j^n).$$
\end{claim}

\proof We first show that there exists a $\chi'(xy,\bar z)\in L$ and $m^*<\omega$ such that
$\forall \bar z \exists ^{<\infty}xy\ \chi'(xy,\bar z)$ and such that for every Morley sequence
$(\bar c_i\vert i\leq m^*)$ of $tp(\bar c/ab)$ that starts at $\bar c$ we have $\chi'(ab, \bar
c_0,...,\bar c_{m^*})$. Indeed, otherwise by Claim \ref{subclaim1} and compactness there is a
Morley sequence $(\bar c^*_i\vert i<\omega)$ of $tp(\bar c/ab)$ that starts at $\bar c$ such that
$ab\not\in acl(\bar c^*_i\vert i<\omega)$. Note that if $e=Cb(\bar c/ab)$ then by our assumption
$e$ is interbounded with $ab$ ($e\in bdd(ab)$ and the assumption $\fork{ab}{\bar c}{bdd(ab)\cap
bdd(\bar c)}$ implies that $\fork{e}{\bar c}{bdd(\bar c)\cap bdd(ab)}$; thus $SU(e/bdd(\bar c)\cap
bdd(ab))>1$ and in particular $SU(e)=2$.) So we get a contradiction to the fact that $e\in dcl(\bar
c^*_i\vert i<\omega)$.

Now, assume by a way of contradiction that the claim is false, i.e. for all $m^*<\omega$ and all
$\chi^*(xy,\bar z)\in L$ (where $\bar z$ has the sort of a tuple of length $m^*+1$ of realizations
of $\bar c$) with $\forall \bar z \exists ^{<\infty}xy\ \chi^*(xy,\bar z)$ there exists a Morley
sequence $(\bar c_i\vert i\leq m^*)$ of $tp(\bar c/ab)$ that starts at $\bar c$ such that
$$\bigwedge_{0\leq i<j\leq m^*} \neg\chi^*(ab,\bar c_0^{<n},\bar c_1^{<n},..., \bar c_{m^*}^{<n},
\bar c_i^n,\bar c_j^n).$$

\noindent By Claim \ref{subclaim1} and compactness, there exists a Morley sequence $(\bar d_i\vert
i\leq m^*)$ of $tp(\bar c/ab)$ that starts at $\bar c$ such that
$$(*)\ \ \ \ \bigwedge_{0\leq i<j\leq m^*} ab\not\in acl(\bar d_0^{<n},\bar d_1^{<n},..., \bar d_{m^*}^{<n},
\bar d_i^n,\bar d_j^n).$$

\noindent By what we saw at the beginning of the current claim,  $$(**)\ \ \ \ \ acl(\bar
d_0^{<n},\bar d_1^{<n},..., \bar d_{m^*}^{<n},ab)=acl(\bar d_0,\bar d_1,..., \bar
d_{m^*},ab)=acl(\bar d_0,\bar d_1,..., \bar d_{m^*}).$$

\noindent Now, let $C\subseteq\bigcup_{0\leq i\leq m^*} \bar d_i$ be a maximal independent set
containing $\bigcup_{0\leq i\leq m^*} \bar d_i^{<n}$. So, since the dimension of $\bigcup_{0\leq
i\leq m^*} \bar d_i$ in the pregeometry $(p^\CC,acl)$ is exactly $(m^*+1)n+2$, we get that $\vert
C\backslash\bigcup_{0\leq i\leq m^*} \bar d_i^{<n}\vert =2$. A contradiction to $(*,**)$ $\Box$.\\

\noindent Let $S=\{(a',b')\in\tilde D\vert\ \mbox{for\ every\ Morley\ sequence} (\bar c_i\vert
i\leq m^*)\ \mbox{of\ } tp(\bar c/a'b')\ \mbox{that\ starts\ at\ } \bar c\ $\\
$$\bigvee_{0\leq i<j\leq m^*} \chi^*(a'b',\bar c_0^{<n},\bar c_1^{<n},...,
\bar c_{m^*}^{<n}, c_i^n,c_j^n)\}.$$

\begin{claim}\label{subclaim3}
Let $S_1$ be the projection of $S$ on the first coordinate. Then $S_1\cup acl(\bar c)$ is an
unbounded Stone-open set over $\bar c$ and for every $a'\in S_1$ there exists an independent tuple
$\bar c'$ of realizations of $p$ such that $\chi_1^*(a',\bar c')$, where $\chi_1^*(x,\bar z)\equiv
\exists y\chi^*(xy,\bar z)$.
\end{claim}

\proof By Claim \ref{subclaim2}, $(a,b)\in S$, and clearly $a\not\in acl(\bar c)$. We conclude that
$S_1$ is unbounded. Note that $S$ is a forking-open set over $\bar c$ ($\tilde D$ is forking open
over $\bar c$ by Remark \ref{ft_remark} and by Claim \ref{subclaim1}, $S$ is the intersection of
$\tilde D$ with a Stone-open set over $\bar c$; therefore $S$ is forking-open over $\bar c$ by
Remark \ref{ft_remark}). Now, by Remark \ref{remark1}, $D_*$ has PCFT and hence $S_1$ is a
forking-open set over $\bar c$. Now, as $SU(S_1)=1$, Fact \ref{tau bounded SU} implies that
$S_1\cup acl(\bar c)$ is a Stone-open set over $\bar c$. Assume now that $a'\in S_1$. Then for some
$b'$, $(a',b')\in S$. So, $\chi^*(a'b',\bar c_0^{<n},\bar c_1^{<n},..., \bar c_{m^*}^{<n},
c_i^n,c_j^n)\}$ for some (all) Morley sequence $(\bar c_k\vert k\leq m^*)$ of $tp(\bar c/a'b')$
that starts at $\bar c$ and some $0\leq i<j\leq m^*$. As $\{a'b',\bar c_0^{<n},\bar c_1^{<n},...,
\bar c_{m^*}^{<n}\}$ is $\emptyset$-independent we conclude, by counting dimensions in the
pregeometry $(D,acl)$, that $\{\bar c_0^{<n},\bar c_1^{<n},..., \bar c_{m^*}^{<n}, c_i^n,c_j^n)\}$
is $\emptyset$-independent and clearly $\chi^*_1(a', \bar c_0^{<n},\bar c_1^{<n},..., \bar
c_{m^*}^{<n}, c_i^n,c_j^n)$.$\ \ \ \ \Box$\\

Now, let $\theta(x)\in L(\bar c)$ be any formula such that $a\models\theta(x)$ and
$\theta^\CC\subseteq S_1\cup acl(\bar c)$. Then $\theta(x)$ is the required formula.$\ \ \ \ \Box$\\

\begin{remark}\em
Note that the set $\tilde D$ defined in the proof of Theorem $\ref{main}$ is not definable unless
the partial type $\Gamma(x,\bar c)=D(x)\wedge [x\not\in acl(\bar c)]$ is implied by consistent
formula over $\bar c$ (e.g. when $D=D(x)$ is strongly minimal, this holds iff $D^\CC\cap acl_x(\bar
c)$ is finite). This is the only reason we are forced to deal with the forking topology in this
proof.
\end{remark}

\proof If $\tilde D$ is definable then $\exists x'\ \tilde D(x,x')$ implies $\Gamma(x,\bar c)$.\qed

As a special case we get Buechler's dichotomy for minimal $D$-rank 1 types in any simple theory.

\begin{corollary}
Let $p\in S(\emptyset)$ be a minimal type with $D(p)=1$. Then either $p$ is 1-based or there exists
$\tilde \chi(x,\bar z)\in L$ with $\forall \bar z\exists ^{<\infty}x \tilde\chi(x,\bar z)$ and
$\theta^*(x)\in p(x)$ such that for any non-algebraic realization (over $\emptyset$) $a$  of $\theta^*(x)$ there is an
$\emptyset$-independent tuple $\bar c$ of realizations of $p$ such that $\tilde\chi(a,\bar c)$. In
particular, $RM(\theta^*(x))=1$.
\end{corollary}

\proof
Assume $p$ is not 1-based. Let $\bar c$, $\theta(x)=\theta(x,\bar c), \tilde\chi(x,\bar z)\in L$ be
given by Theorem \ref{main}; so $\theta(x)\in \bar p$, where $\bar p\in S(\bar c)$ is the
unique non-algebraic complete extension of $p$ over $\bar c$. Now, $p^\CC\subseteq acl(\bar c)\cup
\theta^\CC$. By compactness, there exists $\theta^*(x)\in p$ with $(\theta^*)^\CC\subseteq acl(\bar
c)\cup \theta^\CC$. We may clearly assume that $D(\theta^*(x))=1$. It follows that for every non-algebraic realization  (over $\emptyset$) $a$ of $\theta^*(x)$ there
is an $\emptyset$-independent tuple $\bar c'$ of realizations of $p$ such that $\tilde\chi(a,\bar
c')$ (as  any non-algebraic realization $a$ of $\theta^*(x)$ realize some $\emptyset$-conjugate $\theta(x,\bar c')$ of $\theta(x,\bar c)$).
To see that this implies $RM(\theta^*)=1$, note that for any set $A$, the formula
$\theta^*(x)$ has only finitely many complete non-algebraic extensions over $A$. Indeed, otherwise, there are $a_n, \bar c'_n$ for $n<\omega$
such that $\{tp(a_n/A) \vert  n<\omega\}$ are non-algebraic distinct types that implies $\theta^*(x)$ and such that $\tilde \chi(a_n,\bar
c'_n)$ holds and each $\bar c'_n$ is an $\emptyset$-independent tuple of realizations of $p$.
We may assume that $\bar c'_n$ is independent from $A$  over $a_n$, and in particular, since $D(\theta^*(x))=1$, $\bar c'_n$ is independent from $A$  over $\emptyset$ for every $n<\omega$.
This is a contradiction to the fact that for any $n<\omega$ and any set $B$ there is a unique type over $B$ that doesn't fork over $\emptyset$
whose restriction to $\emptyset$ is the type of $n$ $\emptyset$-independent realizations of $p$ .$\ \ \ \ \Box$\\

We aim now to prove that a $D$-rank 1 type that is topologically close to a 1-based type is in fact
1-based.

\begin{lemma}\label {lemma2}
Let $D$ be a weakly minimal definable set over $\emptyset$. Let $\bar c\subseteq D$ be any tuple.
Then the set $D^2_{\cal NO}(\bar c)=\{(a,b)\in D^2 \vert \fork{ab}{\bar c}{bdd(ab)\cap bdd(\bar c)}
\}$ is a forking-open set over $\bar c$.
\end{lemma}

\proof Assume $(a,b)\in D^2_{\cal NO}(\bar c)$. It will be sufficient to show that there exists a
forking-open set $\UU$ over $\bar c$ such that $(a,b)\in \UU\subseteq D^2_{\cal NO}(\bar c)$. By
Claim \ref{claim1}, and the fact that $D_*$ is supersimple (and thus eliminates hyperimaginaries),
we may work in $D_*^{eq}$ and replace $bdd$ by $acl=acl^{eq}$ in the definition of $D^2_{\cal
NO}(\bar c)$ (note that every forking-open set over some set $A$ in $D_*^{eq}$ is a forking open
set over $A$ in $\CC^{eq}$.) \textbf{So, from now on we work in $D_*^{eq}$}. Clearly we may assume
that $\bar c=c_0c_1...c_n$ is an $\emptyset$-independent. As $(D,acl)$ is a pregeometry, we may
assume, as in the proof of Theorem \ref{main}, that there are $\chi_0(x,y,\bar z)\in L$,
$\chi_1(z_n,z_{n-1},...z_0,x,y)\in L$ such that $\nonforkempty{a}{\bar c}$, $\chi_0(b,a,\bar c)$,
$\forall y\forall \bar z\exists ^{<\infty}x\chi_0(x,y,\bar z)$ and
$\chi_1(c_n,c_{n-1},...,c_0,a,b)$, and
$$\forall z_0z_1...z_{n-1}\forall xy \exists ^{<\infty}z_n\chi_1(z_n,z_{n-1},...z_0,x,y).$$
Let $$\tilde D=\{(a',b')\in D^2\vert\ a'\not\in acl(\bar c), \chi_0(b',a',\bar c),
\chi_1(c_n,c_{n-1},...,c_0,a',b')\}.$$ As before, $\tilde D$ is a forking-open set over $\bar c$,
and $(a,b)\in\tilde D$. Now, since $\fork{ab}{\bar c}{acl(ab)\cap acl(\bar c)}$, there exists $e\in
dcl(Cb(\bar c/ab))\backslash acl(\bar c)$. Let $s$ be the sort of $e$. Let $$\UU=\{(a',b')\in
\tilde D \vert \exists e'\in (D_*^{eq})^s [e'\not\in acl(\bar c)\wedge e'\in acl(Cb(\bar
c/a'b'))]\}.$$

\noindent To finish the proof it remains to show the following.

\begin{claim}\label{subclaim4}
$(a,b)\in \UU$ and $\UU\subseteq D^2_{\cal NO}(\bar c)$, and $\UU$ is a forking-open set over $\bar
c$.
\end{claim}

\proof By the definitions of $e$ and $\UU$, $(a,b)\in \UU$. By the definitions $D^2_{\cal NO}(\bar
c)$ and $\UU$, clearly $\UU\subseteq D^2_{\cal NO}(\bar c)$. To prove that $\UU$ is a forking open
set first note that for every $(a',b')\in \tilde D$, the dimension of $\{c_n,...c_1,c_0,a',b'\}$ in
the pregeometry $(D,acl)$ is $n+2$ , and therefore $\{c_{n-1},...c_1,c_0,a',b'\}$ is
$\emptyset$-independent. In particular, $a'b'$ is independent from $\{c_{n-1},...c_1,c_0\}$ for all
$(a',b')\in\tilde D$. By Remark \ref{Open_Cb_remark}, we conclude that the set $\{(e',a'b') \vert\
e'\in acl(Cb(\bar c/a'b'))\}$ is relatively Stone-open over $\bar c$ in $F=\{(e',a'b') \vert\ e'\in
(D_*^{eq})^s, (a',b')\in\tilde D\}$. By Remark \ref{remark1}, $D_*^{eq}$ has PCFT, thus $\UU$ is a
forking-open set over $\bar c$.$\ \ \ \ \Box$

\qedend

The following notion (from [S2]) is a strong version of the notion of "essential 1-basedness" (the
only difference here is that $\UU$ is not required to be type-definable).


\begin{definition}\label {def ess 1-based}\em
1) A type $p\in S(A)$ is said to be \em s-essentially 1-based (by means of the forking topology)
\em if for every finite tuple $\bar c$ from $p$ and for every forking-open set $\UU$ over $A\bar
c$, the set $\{a\in \UU \vert\ Cb(a/A\bar c)\not\in bdd(aA)\}$ is nowhere dense in the
Stone-topology of $\UU$ \footnote{Note that if $\UU=\emptyset$, the requirement on $\UU$ is
trivially satisfied.}.

\noindent 2) We say that a unidimensional simple theory $T$ is \em s-essentially 1-based \em if for
every $SU$-rank 1 partial type $p_0$ over some $A$, every $p\in S(A)$ that is internal in $p_0$ is
s-essentially 1-based by means of the forking topology.

\end{definition}

\begin{definition}\em
A type $p\in S(A)$ is said to be \em densely 1-based \em if for every finite tuple $\bar c$ of
realizations of $p$ and every forking-open set $\UU$ over $A\bar c$ such that $(p^n)^\CC\cap\UU\neq
\emptyset$ for some $n<\omega$, there exist $a\in\UU$ such that $\nonfork{a}{\bar c}{bdd(A\bar
c)\cap bdd(Aa)}$.
\end{definition}

\begin{remark}\em
If $p\in S(A)$ is s-essentially 1-based then $p$ is densely 1-based.
\end{remark}

\proof By the assumption, if $\UU$ is a forking-open set over $A\bar c$, where $\bar c$ is some
finite tuple of realizations of $p$, then the set $\{a\in \UU \vert\ Cb(a/A\bar c)\not\in
bdd(aA)\}$ is nowhere dense in the Stone-topology of $\UU$. Now, to show that $p$ is densely
1-based, assume $\UU$ is a forking-open set over $A\bar c$ such that $(p^n)^\CC\cap\UU\neq
\emptyset$. In particular, $\UU\not=\emptyset$ and thus the set $\{a\in \UU \vert\ Cb(a/A\bar c)\in
bdd(aA)\}$ is not empty, so there
exist $a\in\UU$ such that $\nonfork{a}{\bar c}{bdd(A\bar c)\cap bdd(Aa)}$.$\ \ \ \ \Box$\\

\begin{theorem}\label{thm2}
Let $p\in S(\emptyset)$ be a type of $D$-rank 1. If $p$ is densely 1-based then $p$ is 1-based. In
particular, if $p\in S(\emptyset)$ is an s-essentially 1-based type of $D$-rank 1 then $p$ is
1-based.
\end{theorem}

\proof Assume $p\in S(\emptyset)$ is a type of $D$-rank 1 and $p$ is densely 1-based. If $p$ is not
1-based then by Lemma \ref{lemma1}, there are $a,b\in p^\CC$ and a finite tuple $\bar c$ of
realizations of $p$ such that $\fork{ab}{\bar c}{bdd(ab)\cap bdd(\bar c)}$. Let $D\in p$ be of
$D$-rank 1. By Lemma \ref{lemma2}, we conclude that $D^2_{\cal NO}(\bar c)$ is a non-empty
forking-open set over $\bar c$. Contradiction to the assumption that $p$ is densely 1-based.$\ \ \ \ \Box$\\

\begin{corollary}
A supersimple unidimensional theory that is s-essentially 1-based is 1-based. In particular, any
countable hypersimple unidimensional theory that is s-essentially 1-based is 1-based.
\end{corollary}

\proof First, recall the following result [W1].

\begin{fact}\label{wagner_fact}
Let $T$ be any simple theory and work with hyperimaginaries. Assume $p\in S(A)$ is analyzable in an
$A$-invariant family of 1-based types. Then $p$ is 1-based.
\end{fact}

\noindent Let $T$ be a unidimensional supersimple theory that is s-essentially 1-based. Let $D$ be
a weakly minimal set (a non-algebraic definable set of minimal $D$-rank). Then any non-algebraic
completion $p$ of $D$ is a type of $D$-rank 1. By Theorem \ref{thm2}, $p$ is 1-based.

By unidimensionality, every type over the domain of $p$ is analyzable in $p$ and thus 1-based by
Fact \ref{wagner_fact}. Thus $T$ is 1-based (note that if $tp(a/A)$ is 1-based and $a$ is
independent from $A$ then $tp(a)$ is 1-based (see [W1]), thus every type over the empty set is
1-based). The last statement follows by supersimplicity of countable
hypersimple unidimensional theories [S1].$\ \ \ \ \Box$\\

\section{Stable types in hypersimple unidimensional theories}
In this section we observe that the existence of a non-algebraic stable partial type in a
hypersimple unidimensional theory implies superstability. $T$ will denote an arbitrary complete theory.\\

The following definition is standard.

\begin{definition}\em
Let $p(x)$ be a partial type over $\CC$.

\noindent 1) $p(x)$ is called \em stable for $\phi(x,y)$ \em if there does not exists a sequence
$(a_ib_i \vert i<\omega)$ such that $a_i\models p$ and such that $\phi(a_i,b_j)$ iff $i<j$.

\noindent 2) $p(x)$ is said to be \em stable \em if $p(x)$ is stable for any formula $\phi(x,y)\in
L$.
\end{definition}

\noindent In this section we show the following.

\begin{proposition}\label{uni_stable}
Let $T$ be a hypersimple unidimensional theory. Assume there exists a non-algebraic stable partial
type. Then $T$ is superstable.
\end{proposition}

The following remark follows easily from well known results.

\begin{remark}\label{remark2}\em
Let $p(x)$ be a partial type. Then the following are equivalent.

\noindent 1) $p=p(x)$ is stable.\\
\noindent 2) For every infinite cardinal $\lambda$ such that $\lambda^{\vert T\vert}=\lambda$, for
every set $A\supseteq dom(p)$ with $\vert A\vert=\lambda$, we have $\vert \{ q\in S(A) \vert\
p\subseteq q\}\vert \leq \lambda$.

\noindent 3) For some infinite cardinal $\lambda\geq \vert dom p\vert$, for every set $A\supseteq
dom(p)$ with $\vert A\vert=\lambda$, $\vert \{ q\in S(A) \vert\ p\subseteq q\}\vert \leq \lambda$.

\noindent 4) For every model $M$, every type $q\in S(M)$ that extends $p$ is definable.

\noindent 5) Every type $q\in S(B)$ that extends $p$ is definable.

\noindent 6) For every formula $\phi=\phi(x,y)$, $R(p,\phi,2)<\omega$.

\end{remark}

\proof $1)\Rightarrow (4$ follows by the usual proof of definability of $\phi$-types over a model,
using the fact that if $p(x)$ is stable for $\phi(x,y)$ then there is a $\psi(x)\in p(x)$ such that
$\psi(x)$ is stable for $\phi(x,y)$ (compactness). $4)\Rightarrow (2$ is clear. $2)\Rightarrow (3$
is trivial. $3)\Rightarrow (6$: otherwise let $\lambda$ be as given in 3) and let $\mu$ be the
minimal cardinal such that $2^\mu>\lambda$ and let $\{a_\eta \vert \eta\in 2^{<\mu}\}$ be such that
for every $\bar\eta\in 2^\mu$, $p(x)\wedge\bigwedge_{i<\mu}{\phi(x,a_{\bar\eta\vert
i}})^{\bar\eta(i)}$ is consistent, contradicting the assumption in 3). $6)\Rightarrow (5$: given
$q(x)\in S(B)$ extending $p(x)$ and any $\phi=\phi(x,y)\in L$ let $r=R(q(x),\phi,2)$ and let
$\psi(x)\in q(x)$ such that $r=R(\psi(x),\phi,2)$. Then, for any $b\in B$ we have $\phi(x,b)\in q$
iff $R(\psi(x)\wedge\phi(x,b),\phi,2)=r$. This is a definable condition on $b$ over $B$.
$5)\Rightarrow (2$ is clear. $6)\Rightarrow (1$ Otherwise there exists a sequence $(a_ib_i \vert
i\in \mathbb{Q})$ such that $a_i\models p$ and such that $\phi(a_i,b_j)$ iff $i<j$. The required
consistency is now obvious.$\ \ \ \ \Box$\\

\begin{remark}\label{remark3}\em
\noindent 1) Assume $p_i(x_i)$ for $i<n$ ($n<\omega$) are stable partial types over $\CC$. Then so
is $\bigwedge_{i<n} p_i(x_i)$.

\noindent 2) Assume $p(x),q(y)$ are partial types over $\CC$ such that for some small set $B$ we
have $q^\CC\subseteq acl(p^\CC\cup B)$. Then, if $p(x)$ is stable so is $q(y)$.

\noindent 3) Assume $\Gamma(x)\equiv\bigvee_{i} p_i(x)$, where each of $\Gamma(x), p_i(x)$ are
partial types over $\CC$, and assume $p_i(x)$ are stable. Then $\Gamma(x)$ is stable.
\end{remark}

\proof 1) We may clearly assume $n=2$. In this case let $B$ be a set containing the domains of both
$p_0(x_0)$ and $p_1(x_1)$, so the type of $(x_0,x_1)$ over $B$ is determined by the type of $x_0$
over $B$ and of $x_1$ over $Bx_0$. By stability of $p_0(x_0)$ and $p_1(x_1)$ we are done.

\noindent 2) Let $A$ be a sufficiently large superset of $B$ and of the domains of $p,q$ with
$\vert A\vert^{\vert T\vert}=\vert A\vert$. Then for every $c\in q^\CC$ there is an algebraic
formula $\chi(y,x_0,...x_n,b)\in L(B)$ in $y$, and realizations $a_0,...,a_n$ of $p$ such that
$\models \chi(c,a_0,...,a_n,b)$. By 1) and Remark \ref{remark2} the number of possible types of the
tuple $(a_0,...a_n)$ over $A$ is $\leq\vert A \vert$ and thus so is the number of possible types of
$c$ over $A$.

\noindent 3) is immediate.$\ \ \ \ \Box$\\

\begin{lemma}\label{int_stable}
Let $T$ be simple. Assume $q$ is a stable partial type over $A$, and $p$ is a partial type over $A$
that is almost $q$-internal. Then $p$ is stable.
\end{lemma}

\proof Recall the following basic fact about simple theories [W, Proposition 3.4.9] or [S3, Theorem
5.6]; for a stronger version see also [S4, Theorem 2.2].

\begin{fact}\label{a-gen}
Let $p$ be a partial type over $A$ that is almost-internal in an $A$-invariant set $\UU$ for some
small set $A$. Then there exists a set $B$ such that $p^\CC\subseteq acl(\UU\cup B)$.
\end{fact}

By Fact \ref{a-gen}, there exists a small set $B$ such that $p^\CC\subseteq acl(q^\CC\cup B)$. By
Remark \ref{remark3}(2), we are done.$\ \ \ \ \Box$\\

\begin{lemma}\label{lemma_int_st}
Let $T$ be simple. Assume $q=q(x,b)$ is a stable partial type. Let $A\supseteq b$ be a small
infinite set such that $\vert A\vert^{\vert T\vert+\vert b\vert}=\vert A\vert$. Then  $$\vert\{p\in
S(A)\vert\ p\ is\ almost\ q-\mbox{internal\ }\}\vert\leq \vert A\vert.$$
\end{lemma}

\proof Assume $p\in S(A)$ is almost $q$-internal. Let $A_0\subseteq A$ be such that $p$ doesn't
fork over $A_0$ and with $b\subseteq A_0$, $\vert A_0\vert\leq \vert T\vert+\vert b\vert$. Now,
$p_0=p\vert A_0$ is almost $q$-internal and thus by Lemma \ref{int_stable}, $p_0$ is stable. In
particular, every $p\in S(A)$ that is almost $q$-internal extends a stable type over a subset of
$A$ of size $\leq \vert T\vert+\vert b\vert$. By Remark \ref{remark2}(2) and the fact that the
number of types over subsets of $A$ of size $\leq \vert T\vert+\vert b\vert$ is $\leq\vert A\vert$,
we are done.$\ \ \ \ \Box$\\

\begin{corollary}\label{st_cor}
Let $T$ be hypersimple. Assume $q=q(x,b)$ is a stable partial type and let $A\supseteq b$ be a
small set.

\noindent 1) If $A$ is infinite such that $\vert A\vert^{\vert T\vert+\vert b\vert}=\vert A\vert$,
then
$$\vert\{p\in S(A)\vert\ p\ is\ q-\mbox{analyzable\ }\}\vert\leq \vert A\vert.$$
\noindent 2) Assume $p\in S(A)$ is analyzable in $q$ (by an imaginary sequence). Then $p$ is
stable.
\end{corollary}

\proof To prove 1), it will be sufficient to note that for every $\alpha<{\vert T\vert}^+$
$$\vert \{tp((a_i \vert i\leq\alpha)/A) \vert\ (a_i\vert i\leq\alpha)\ \mbox{is\ an\ analysis\ in}\ q\ \mbox{over}\
A\}\vert\leq \vert A\vert $$ (we say that $(a_i\vert i\leq\alpha)$ is an analysis in $q$ over $A$
if $tp(a_i/(a_i\vert j<i)\cup A)$ is  $q$-internal for all $i\leq \alpha$). Indeed, this follows by
repeated applications of Lemma \ref{lemma_int_st} and the assumption that $A$ is infinite and
$\vert A\vert^{\vert T\vert+\vert b\vert}=\vert A\vert$ (which implies $\vert A\vert^{\vert
T\vert}=\vert A\vert$). Now, as in particular $\vert A\vert\geq {\vert T\vert}^+$, we get the
required statement.

To prove 2), assume $p\in S(A)$ is analyzable in $q$ and let $B\supseteq A$ be any infinite small
set such that $\vert B\vert^{\vert T\vert+\vert b\vert}=\vert B\vert$. It will be sufficient to
show that the number of complete extensions of $p$ over $B$ is $\leq \vert B\vert$. Indeed, let $a$
be any realization of $p$. Then there is an analysis $(a_i\vert i\leq \alpha)$ in $q$ over $A$ for
some $\alpha<{\vert T\vert}^+$ such that $a_\alpha=a$. Therefore $(a_i\vert i\leq\alpha)$ is an
analysis in $q$ over $B$. Thus, $tp(a/B)$ is analyzable in $q$ over $B$. By part 1) we conclude
that the number of extensions of $p$ over $B$ is $\leq\vert B\vert$, thus $p$ is stable.$\ \ \ \ \Box$\\

\noindent\textbf{Proof of Proposition \ref{uni_stable}:} Let $q(x,b)$ be a stable non-algebraic
partial type. By the assumption that $T$ is a hypersimple unidimensional theory, every complete
type over a superset of $b$ is analyzable in $q$. By Corollary \ref{st_cor}(2), every complete type
over a superset of $b$ is stable. Thus $T$ is stable over $b$ and so $T$ is stable. By [H], $T$ is
superstable.$\ \ \ \ \Box$\\

Here is an application of Proposition \ref{uni_stable}.

\begin{corollary}
Let $T$ be a hypersimple non-stable unidimensional theory. Then $T$ is non-trivial. In particular,
every countable trivial hypersimple unidimensional theory is $\aleph_1$-categorical.
\end{corollary}

\proof Assume $T$ is a hypersimple unidimensional trivial theory. Let $p\in S(A)$ be a $SU$-rank 1
type with $\vert A\vert=\vert T\vert$. By Proposition \ref{uni_stable} it will be sufficient to
prove that $p$ is stable. Indeed, by Remark \ref{remark2}, we need to show that for every
$B\supseteq A$ of size $\vert T\vert$, the number of extensions of $p$ to $B$ is $\leq \vert
T\vert$. To see that, let $p_0,p_1\in S(B)$ be two non-algebraic extensions of $p$. Since $p_0,p_1$
are non-orthogonal there is a set $C\supseteq B$ and $a_0\models p_0, a_1\models p_1$ such that
each $a_i$ is independent from $C$ over $B$ and such that $a_0$ depends on $a_1$ over $C$. By
triviality of $T$, it follows that $a_0$ depends on $a_1$ over $B$. Thus we have shown that any two
non-algebraic extensions of $p$ over $B$ have interalgebraic realizations over $B$ so the number of
extensions of $p$ over $B$ is bounded by $\vert B\vert=\vert T\vert$.

\noindent Now, if $T$ is a countable trivial hypersimple unidimensional theory then by the first
part $T$ must be stable. Let $p^*$ be a minimal type. Since $T$ is trivial, no infinite group is
interpretable and thus $\CC=acl({p^*}^\CC)$ (every element in $\CC$ is analyzable in $p^*$, but
since no binding group is infinite we get algebraicity in each step). Since $L$ is countable, $T$
is $\omega$-stable and thus $\aleph_1$-categorical.

\noindent Ziv Shami, E-mail address: zivshami@gmail.com\\
Dept. of Mathematics and Computer Science\\
Ariel University\\
Samaria, Ariel 44873\\
Israel.\\

\end{document}